\documentclass[11pt]{article}
\usepackage{amssymb,amsfonts,amsthm}
\usepackage{color}
\usepackage[latin1]{inputenc}
\usepackage[active]{srcltx}
\def\openC{{\rm C\kern-.18cm\vrule width.8pt height 7pt depth-.2pt \kern.18cm}}
\def\openN{{{\rm I}\kern-.16em {\rm N}}}
\def\openR{{{\rm I}\kern-.16em {\rm R}}}
\def\openT{{{\rm T}\kern-.42em {\rm T}}}
\def\openZ{{{\rm Z}\kern-.28em{\rm Z}}}

\def\eop{\hfill\rule{2.5mm}{2.5mm}}
\def\pf{\par\smallbreak\noindent {\bf Proof.} \ }

\newtheorem{thm}{Theorem}[section]

\newtheorem{lem}[thm]{Lemma}
\newtheorem{prop}[thm]{Proposition}

\theoremstyle{definition}

\newtheorem{rem}[thm]{Remark}
\newtheorem{ex}[thm]{Example}

\def\eop{\hfill\rule{2.5mm}{2.5mm}}
\begin{document}

\title{
{\textbf{\Large{Kolmogorov widths on the sphere via eigenvalue estimates for H\"{o}lderian integral operators
%Approximating operators: their application in Kolmogorov $n$-widths via eigenvalues problems on the spherical setting
}}} \vspace{-4pt}
\author{
T. Jord\~{a}o\,\thanks{
Both authors partially supported by FAPESP, grants $\#$ 2016/02847-9.\ First author also supported by the additional grant 2017/07442-0. }
\,\,\&\,\,
V. A. Menegatto}
}

\date{}
\maketitle \vspace{-30pt}
\bigskip

\begin{center}
\parbox{13 cm}{{\small Approximation processes in the reproducing kernel Hilbert space associated to a continuous kernel on the unit sphere $S^m$ in the Euclidean
space $\mathbb{R}^{m+1}$ are known to depend upon the Mercer's expansion of the compact and self-adjoint $L^2(S^m)$-operator associated to the kernel.\ The estimation of the Kolmogorov $n$-th width
of the unit ball of the reproducing kernel Hilbert space in $L^2(S^m)$ and the identification of the so-called optimal subspace usually suffice.\ These Kolmogorov widths can be computed through the eigenvalues of the integral operator associated to the kernel.\ This paper provides sharp upper bounds for the Kolmogorov widths in the case in which the kernel satisfies an abstract H\"{o}lder condition.\ In particular, we follow the opposite direction usually considered in the literature, that is, we estimate the widths from decay rates for the sequence of eigenvalues of the integral operator.}}
\end{center}

%\noindent{\bf Keywords:} convolutions, decay rates, eigenvalues, integral operators,
%singular values, spectral analysis.\\
%{\bf MSC:} 41A36, 42A82, 45C05, 45M05, 45P05, 47A75, 47B34, 47B38, 47G10.

%\noindent{\bf Keywords: convolution, Fourier coefficients, integral operator, multiplier operator}

\thispagestyle{empty}

%
%-------------------------------------------------------------------------------------------------------------------------------------------------------
%

\section{Introduction}\label{s1}

Let us start with some background material.\ We will endow the unit sphere $S^m$, $m\geq2$, of $\mathbb{R}^{m+1}$ with its usual geodesic distance and write $\sigma_m$ to denote the induced Lebesgue measure on $S^m$ and $\omega_m$ its volume.\ If $K: S^m \times S^m \to \mathbb{R}$ is a symmetric and positive definite kernel on $S^m$, write $(\mathcal{H}(K), \| \cdot \|_{\mathcal{H}})$ to denote the unique separable Hilbert space of functions $f: S^m \to \mathbb{R}$ where $K$ is a reproducing kernel.\ If $K$ is continuous, the space $\mathcal{H}(K)$ is embeddable in the usual space $L^2(S^m):=L^2(S^m,\sigma_m)$.\ Indeed, this follows from
$$
\int_{S^m}K(x,x)d\sigma_m(x)<\infty,
$$
and the inequality
$$
\|f\|_2 \leq \left[\frac{1}{\omega_m}\int_{S^m}K(x,x)d\sigma_m(x)\right]^{1/2} \|f\|_{\mathcal{H}}, \quad f\in \mathcal{H}(K).
$$
By the way, we will write $\|\cdot \|_p$, $1\leq p\leq \infty$, to the denote the $p$-norm in the usual space $L^p(S^m):=L^p(S^m,\sigma_m)$.

Under the setting in the previous paragraph, the integral operator $\mathcal{K} : L^2(S^m) \to L^2(S^m)$ given by
\begin{equation}\label{intequ}
\mathcal{K}(f) = \int_{S^m} K(x,y) f(y) d\sigma_m(y), \quad f \in L^2(S^m),
\end{equation}
is well-defined, compact, and self-adjoint.\ Its range is a dense subset of $\mathcal{H}(K)$ and, in addition,
$$
\langle f,g \rangle_{L^2(S^m)}= \langle f,g\rangle_{\mathcal{H}(K)}, \quad f\in \mathcal{H}(K), \quad g \in L^2(S^m).
$$
Since a version of the classical Mercer's Theorem hold, the integral operator $\mathcal{K}$ is positive and has a countable set of positive eigenvalues, say, $\lambda_1 \geq \lambda_2\geq \cdots >0$, with respective eigenfunctions
$\varphi_1, \varphi_2, \ldots$, that is,
$$
\mathcal{K} \varphi_i=\lambda_i \phi_i, \quad i=1,2,\ldots.
$$
The set $\{\varphi_i: i=1,2,\ldots\}$ is orthonormal in $L^2(S^m)$ and orthogonal in $\mathcal{H}(K)$.\ Further,
\begin{equation}\label{mercer}
K(x,y)=\sum_{i=1}^{\infty}\lambda_i\varphi_i(x)\varphi_i(y), \quad x,y\in S^m,
\end{equation}
where the sum is absolutely and uniformly convergent.\ Since $\|\varphi_i\|_{\mathcal{H}(K)}^2=\lambda_i^{-1}$, $i=1,2,\ldots$, it follows that the set
$\{\sqrt{\lambda_i} \, \varphi_i: i=1,2,\ldots\}$ is an orthonormal basis of $\mathcal{H}(K)$.\ Throughout, we will write
$$
\Omega_n:= \mbox{span}\,\{\sqrt{\lambda_i} \, \varphi_i : i=1,\ldots, n\},\quad n=1,2,\ldots.
$$

This is the point at which we may say a little bit about widths of Kolmorogov (see \cite{pinkus} and other references quoted there).\  The {\em Kolmogorov $n$-width} of a subset $A$ of a Hilbert space $H$ is the quantity $d_n(A;H)$ that measures how $n$-dimensional subspaces of $H$ can approximate $A$.\
In other words, it is defined as
\begin{equation}\label{kolmogorovn}
d_n(A;H):= \inf_{V_n\subset H} \sup_{f\in A}  \inf_{f_n\in V_n}  \|f-f_n\|_{H},
\end{equation}
where the first infimum is taken over all $n$-dimensional subspaces $V_n$ of $H$.\ If the infimum is attained, that is,
$$
d_n(A;H)=\sup_{f\in A}  \inf_{f_n\in V_n}  \|f-f_n\|_H
$$
for some $n$-dimensional subspace $V$ of $H$, then $V$ is called an {\em optimal subspace}.\ The characterization of optimal subspaces and either the computation or estimation of the widths are the highlight problems in this regard and, usually, the case in which $A$ is the closed unit ball in $H$ receives most of the attention.

Returning to the spherical setting we previously introduced and letting $S$ be the unit sphere in $\mathcal{H}(K)$, a result in \cite{schaback} (see also \cite[Chapter 6]{pinkus}) reveals that
$$
d_n(S;L^2(S^m))=\inf_{V_n\subset L^2(S^m)}  \sup_{f\in S} \|f-Q_n(f)\|_{2}=\sqrt{\lambda_{n+1}},
$$
in which $Q_n(f)$ is the projector of $f$ onto $V_n$ in $L^2(S^m)$, that is,
$$
Q_n(f)=\sum_{i=1}^n \langle f, h_i\rangle_{L^2(S^m)} h_i, \quad f \in L^2(S^m),
$$
where $\{h_1, h_2, \ldots, h_n\}$ is an $L^2(S^m)$-orthonormal basis of $V_n$.\ In addition, $\Omega_n$ turns out to be the unique optimal subspace.\ In particular, the analysis of the Mercer's representation for $K$ has extreme relevancy in the understanding of the approximation processes in $\mathcal{H}(K)$ which are dictated by the Kolmogorov widths.

If we replace the projector $Q_n$ with the projector $P_n$ of $V_n$ onto $\mathcal{H}(K)$, but keep the approximations in the $L^2(S^m)$ norm, a result in \cite{santin} ratifies that
$$
d_n(S;L^2(S^m))=\inf_{V_n\subset \mathcal{H}(K)}  \sup_{f\in S}  \|f-P_n(f)\|_{2}.
$$
Further, the optimality of $\Omega_n$ remains for this alternative definition of the Kolmogorov $n$-with, that is,
$$
d_n(S;L^2(S^m))= \sup_{f\in S} \,\,\,\, \left\|f-\sum_{i=1}^n \langle f, g_i\rangle_{L^2(S^m)} g_i \right\|_{2}.
$$
where $\{g_1, g_2, \ldots, g_n\}$ is an $\mathcal{H}(K)$-orthonormal basis of $\Omega_n$.

In this paper, we will provide sharp estimates for the Kolmogorov $n$-with $d_n(S;L^2(S^m))$ described above under the assumption that the kernel $K$ satisfies
an abstract H\"{o}lder condition.\ In Section 2, we introduce notation until the point we are able to introduce the H\"{o}lder condition to be used in the paper which is defined through convolutions with parameterized family of measures.\ We provide a few examples frequently used as concrete realizations for the H\"{o}lder condition and also include a Gaussian-like kernel that satisfies one of the exemplified realizations.\ In Section 3, we introduce a family of approximating operators and provide reasonable conditions in order that these operators be uniformly bounded and have finite rank.\ The operators are used to estimate the eigenvalues of the H\"{o}lderian integral operators $\mathcal{K}$ in Section 4, which lead to sharp estimates for $d_n(S;L^2(S^m))$.\ Section 5 contains a concrete case that exemplifies our achievements.

\section{Convolution with measures and the H\"{o}lder condition}

In this section, we introduce more notation and briefly discuss the H\"{o}lder condition we intend to make use of.

Let $\mathcal{M}_{\rho}(S^m)$ be the set of finite regular measures on $S^m$ which are invariant under the group of rotations of $S^m$ fixing a \textit{pole} $\rho$.\ It becomes a Banach space under the norm
$$
|\mu|(S^m)=\sup\left\{\frac{1}{\omega_m}\left|\int_{S^m}fd\mu\right|: f\in L^1(S^m,\mu);\,\|f\|_{1}\leq 1\right\},
$$
where $|\mu|$ is the total variation of $\mu$.\ If $x\in S^m$, $\mathcal{O}_x^{\rho}$ will denote a rotation of $S^m$ such that $\mathcal{O}_x^{\rho}(x)=\rho$.\ That being said,
we define $\varphi_x : \mathcal{M}_{\rho}(S^m) \to \mathcal{M}_{\rho}(S^m)$ by the formula
$$
\varphi_x(\mu):= \mu\circ \mathcal{O}_x^{\rho},  \quad \mu\in \mathcal{M}_{\rho}(S^m).
$$
The notations above agree with those in \cite{berens,dunkl}.

Next we introduce the notion of isotropy for kernels.\ A kernel $K: S^m \times S^m \to \mathbb{R}$ is {\em isotropic} whenever there is a function $K_i: [-1,1]\longrightarrow \mathbb{R}$ (the isotropic part of $K$) so that
$$
K(x,y)=K_i(x\cdot y), \quad x,y\in S^m,
$$
where $\cdot$ stands for the usual inner product of $\mathbb{R}^{m+1}$.\ If $K \in L^1(S^m\times S^m):=L^1(S^m\times S^m, \sigma_m\times \sigma_m)$ is isotropic, then its norm, also denoted by $\|\cdot \|_1$, can be computed through the formula
$$
\|K\|_{1}:=\frac{\omega_{m-1}}{\omega_m}\int_{-1}^{1}|K_i(u)|(1-u^2)^{(m-2)/2}du.
$$

We now recall a result proved in \cite{dunkl}.

\begin{prop} $(1\leq p\leq \infty)$\label{convo}
Let $\rho$ be a pole in $S^m$.\ If $f$ belongs to $L^p(S^m)$ and $\mu$ is an element of $\mathcal{M}_{\rho}(S^m)$, then the formula
\begin{equation}\label{convolution}
(f\ast\mu)(x):=\frac{1}{\sigma_m}\int_{S^m}f(y)d\varphi_x(\mu)(y),
\end{equation}
defines an element of $L^p(S^m)$ satisfying $\|f\ast \mu\|_p\leq \|f\|_p|\mu|$.\ Further, if $f$ is isotropic, then so is $f\ast\mu$.
\end{prop}

We call $f\ast\mu$ {\em the spherical convolution of $f$ and $\mu$}.\ If we consider a family $\{\mu_t: t\in (0,\pi)\}$ in $\mathcal{M}_{\rho}(S^m)$, then we can define a family $\{T_t : t\in(0,\pi)\}$ of linear operators on $L^p(S^m)$, where each $T_t$ is defined through the convolution just introduced:
$$
T_t(f)= f * \mu_t, \quad f \in L^p(S^m).
$$
It is standard to verify that each $T_t$ is bounded.\ Indeed, the total variation of the $\mu_t$ determines an upper bound for the norm of $T_t$ in the sense that
(see, for example, \cite{berens})
$$
\|T_t\| \leq |\mu_t|,\quad t\in(0,\pi).
$$

If for $u \in [-1,1]$ we write $d\lambda_m(u):=(1-u^2)^{(m-2)/2}du$, then we can also construct a family $\{T_t : t\in(0,\pi)\}$ via the natural embedding
$$
f \in L^1([-1,1], \lambda_m) \hookrightarrow \mu_f \in \mathcal{M}_{\rho}(S^m),
$$
 where
$$
d\mu_f(u)=f(u)(1-u^2)^{(m-2)/2}du,\quad u \in [-1,1].
$$
If we start with a family of isotropic kernels $\{K^t: t\in (0,\pi)\}$ in $L^1(S^m)$, since
$K_i^t \in L^1([-1,1],\lambda_m)$, we can now put
$$
T_t(f)=f* \mu_{ K_{i}^{t}}, \quad f \in L^p(S^m).
$$
In this case the norm inequality for the family becomes
$$
\|T_t\| \leq \|K^t\|_{1}, \quad t \in (0,\pi).
$$
The forthcoming results will be formulated based on the two constructions introduced above.

To proceed, let us write $\mathcal{H}_k^m$ to denote the space of all spherical harmonics of degree $k$ in $m+1$ variables and denote its dimension by $d_k^{m}$.\ It is well-known that
\begin{equation}\label{ineqb}
d_k^{m} \leq 2k^m,\quad k\geq k_0,
\end{equation}
where $k_0=k_0(m)$.\
The orthogonal decomposition $L^2(S^m)=\oplus_{k=0}^{\infty}\mathcal{H}_k^m$ is also well-known while the orthogonal projection of $L^2(S^m)$ over a single $\mathcal{H}_k^m$ is given by the formula
\begin{equation}\label{proj}
\mathcal{Y}_{k}(f)(x)=\frac{d_k^{m}}{\omega_m}\int_{S^m} P_k^{(m-1)/2}(x\cdot y)f(y)d\sigma_m(y), \quad f \in L^2(S^m),\quad x\in S^m,
\end{equation}
in which $P_k^{(m-1)/2}$ is the usual Gegenbauer polynomial of degree $k$ associated to the dimension $m$ and normalized as $P_k^{(m-1)/2}(1)=1$.\

The action of the projections over convolutions is given by
$$
{\cal Y}_k(f\ast\mu)(x)=\mu_k{\cal Y}_k(f)(x), \quad f\in L^2(S^m),\quad \mu\in \mathcal{M}_{\rho}(S^m),
$$
where
$$
\mu_k:={\cal Y}_k(\mu) = \frac{1}{\omega_m}\int_{S^m}P_k^{(m-1)/2}(\rho\cdot y)d\mu(y), \quad k=0,1,\ldots.
$$
The sequence $\{\mu_k\}_{k=1}^\infty$ will be called the \textit{multiplier} of $\mu$.\ Additional information on this specific topic can be found in \cite{berens,dunkl}.

At this point it is important to consider a few concrete examples.

\begin{ex} \label{shifting}\textit{(Shifting operator)} The usual shifting operator is defined by the formula (\cite{berens})
$$
S_tf(x)=\frac{1}{R_m(t)}\int_{R_x^t}f(y)d\sigma_r(y),\quad x \in S^m, \quad f\in L^2(S^m),\quad t\in (0,\pi),
$$
in which $d\sigma_r(y)$ is the volume element of the rim $R_x^t:=\{y \in S^m: x\cdot y=\cos t\}$ and $R_m(t)=\omega_{m-1}(\sin t)^{m-1}$ is its total volume.\ Its convolution structure is defined as
$$
S_t(f)=f\ast \mu_t, \quad t\in (0,\pi), \quad f \in L^2(S^m),
$$
where $\{\mu_t: t\in (0,\pi)\}\subset \mathcal{M}_{\rho}(S^m)$ satisfies
$$
\mathcal{Y}_k(\mu_t)=P_k^{(m-1)/2}(\cos t),\quad k=0,1,\ldots.
$$
In particular, the multiplier of $\mu_t$ is $\{P_k^{(m-1)/2}(\cos t)\}_{k=0}^\infty$.
\end{ex}

\begin{ex}\label{caps}\textit{(Averages on caps)} This example
is discussed in \cite{berens,ditzian2}, while the point of view we will give here is aligned with \cite{jordao1}.\
The average operator on the cap $C_t^x=\{w\in S^m: x \cdot y\geq \cos t\}$
of $S^m$, defined by $t$, is the operator $A_t$ given by
$$
(A_t f)(x)=\frac{1}{C_m(t)}\int_{C_t^x}f(w)d\sigma_m(w), \quad x\in S^m,\quad  t\in (0,\pi),
$$
in which
$C_m(t)$ is total volume of the cap $C_t^x$.\ It is shown in \cite{jordao1} that
$$
A_t(f)=f\ast \mu_{\mathcal{Z}_{t}},\quad t \in (0,\pi), \quad f \in L^2(S^m),
$$
where
$$
\mu_{\mathcal{Z}_{t}}=C_m^{-1}(t)\tilde{\mu}_{\mathcal{Z}_t}, \quad t \in (0,\pi),
$$
$$
d\tilde{\mu}_{\mathcal{Z}_t}(x)=\mathcal{Z}_{t}(\rho,x)d\sigma_m(x), \quad t\in (0,\pi), \quad x \in S^m,
$$
and
$$
\mathcal{Z}_{t}(x, y):=\left\{\begin{array}{rc}
\omega_m,& \quad \mbox{if\ \ } \cos t\leq x\cdot y\leq 1\\
0,& \quad otherwise.
\end{array}\right.
$$
The sequence of projections of $\mu_{\mathcal{Z}_{t}}$ are given by (see \cite{berens})
$$
\mathcal{Y}_k(\mu_{\mathcal{Z}_{t}})=\frac{\omega_{m-1}}{C_m(t)}\left(\int_{0}^{t}P_k^{(m-1)/2}(\cos h)(\sin h)^{m-1}dh\right),\quad t \in (0,\pi), \quad k=0,1,\ldots
$$
while some identities for Gegenbauer polynomials lead to
$$
\mathcal{Y}_k(\mu_{\mathcal{Z}_{t}})=\frac{\omega_{m-1}(m-1)}{k(k+m-1)}(\sin t)^m P_{k-1}^{(m+1)/2}(\cos t),\quad t \in (0,\pi), \quad k=0,1,\ldots,
$$
which defines the multiplier of $\mu_{\mathcal{Z}_{t}}$.
\end{ex}

\begin{ex}\label{stekelov} \textit{(Stekelov-type means)} The Stekelov-type mean defined by $t \in (0,\pi)$ is given by
$$
E_t(f)(x)=\frac{1}{D_m(t)}\int_0^{t}\frac{C_m(s)}{R_m(s)} A_s(f)(x)ds, \quad  t\in (0,\pi), \quad x\in S^m,
$$
where the normalizing constant $D_m(t)$ is chosen so that $E_t(1)=1$.\ We have that
$$
E_t(f)=f \ast \mu_{\mathcal{W}_{t}}, \quad  t \in (0,\pi), \quad f \in L^2(S^m),
$$
where
$$
\mu_{\mathcal{W}_{t}}:=D_m^{-1}(t)\tilde{\mu}_{\mathcal{W}_{t}}, \quad t \in (0,\pi),
$$
$$
d\tilde{\mu}_{\mathcal{W}_{t}}(x)=\mathcal{W}_{t}(\varepsilon\cdot x)d\sigma_m(x), \quad t \in (0,\pi), \quad x \in S^m,
$$
and
$$
\mathcal{W}_{t}(x, y):=\left\{\begin{array}{rc}
\displaystyle{\int_0^t \frac{1}{R_m(s)}\mathcal{Z}_{s}(x,y)ds},& \quad \mbox{if\ \ } \cos t \leq x\cdot y\leq  1\\
0,& \quad otherwise,
\end{array}\right.
$$
$\mathcal{Z}_{s}$ being the kernels described in the previous example.\
Also,
$$
\mathcal{Y}_k(\mu_{\mathcal{W}_{t}})=\frac{1}{D_m(t)}\int_0^{t}\frac{C_m(s)\mathcal{Y}_k(\mu_{\mathcal{Z}_{s}})}{R_m(s)}ds, \quad t\in (0,\pi)\quad k=0,1,\ldots.
$$
and the multiplier of $\mu_{\mathcal{W}_{t}}$ becomes clear.
\end{ex}

Finally, we may introduce the H\"{o}lder condition we intend to use in the paper.\ It depends upon a fixed sequence of measures $\{\mu_t\}$ defining convolution operators $T_t$ as previously described and also exemplified above.\
It also depends upon a real number $\rho\in (0,2]$ and a function $B: S^m \to [0,\infty)$ belonging to $L^\infty(S^m)$.\ A kernel $K: S^m\times S^m\longrightarrow \mathbb{R}$ is said to be ($\mu_t, B, \rho$)-\textit{H\"{o}lder} if
\begin{equation}\label{Holdercond}
|(K(x,\cdot) *\mu_t)(y) - K(x,y)|\leq B(x)t^{\rho}, \quad t\in(0,\pi), \quad x,y\in S^m.
\end{equation}
In this case the corresponding integral operator $\mathcal{K}$ is called {\em H\"{o}lderian}.

The three examples described above define potential sequences of measures that can be used in the H\"{o}lder condition just defined.\ It is easy to verify that if a kernel $K$ is $(\mu_t,B,\rho)$-H\"{o}lder, for some $B$ and $\rho$ according to Example \ref{shifting}, then it is $(\mu_{\mathcal{Z}_{t}},B,\rho)$-H\"{o}lder.\ Similarly, if it is $(\mu_{\mathcal{Z}_{t}},B,\rho)$-H\"{o}lder, then it is $(\mu_{\mathcal{W}_{t}},B,\rho)$-H\"{o}lder as well.

We close the section presenting a concrete realization for the previous definition.\ It is commonly used in learning theory and in methods related to the approximation of functions in reproducing kernel Hilbert spaces (see \cite{cucker, minh} and references quoted there).\ As usual, we will write $A(t)\asymp B(t)$, $t \in (0,\infty)$, to indicate that there exist nonnegative constants $c_1$ and $c_2$, not depending upon $t$, such that $c_1\,A(t)\leq B(t)\leq c_2\, A(t)$, $t \in (0,\pi)$.

\begin{ex} \label{gaussianex} For $\sigma>0$, let $K_{\sigma}$ be the Gaussian-like kernel given by
$$
K_{\sigma}(x,y)=\exp(-2\sigma^{-2}(1-x \cdot y)), \quad x,y \in S^m.
$$
Proposition 2.14 in \cite{cucker} ratifies that $K$ is representable in the form
$$
K_{\sigma}(x,y)=\sum_{k=0}^{\infty}\lambda_k^{\sigma}\sum_{j=1}^{d_k^m}\varphi_{k,j}(x)\varphi_{k,j}(y),
\quad x,y \in S^m,
$$
where $\{\lambda_k^{\sigma}\}_{k=0}^\infty \subset [0,\infty)$ and $\{\varphi_{k,j}: j=1,2\ldots, d_k^m \}$ is an $L^2(S^m)$-orthonormal basis of $\mathcal{H}_k^m$.\ As a matter of fact, we have that
$$
\lambda_k^{\sigma}=e^{-2/{\sigma^2}} \sigma^{m-1}I_{k+(m-1)/2}(2\sigma^{-2})\Gamma((m+1)/2), \quad k=0,1,\ldots.
$$
where $I_v(\cdot)$ stands for the modified Bessel function of first kind associated with $v$.\ In particular,
$$
2\lambda_{k}^{\sigma}>(2k+m+1)\sigma^2\lambda_{k+1}^{\sigma}, \quad
k=0,1,\ldots.
$$
Since
\begin{eqnarray*}
S_t(K_{\sigma}(x, \cdot))=\sum_{k=0}^{\infty}\lambda_k^{\sigma}P_k^{(m-1)/2}(\cos t)\sum_{j=1}^{d_k^m}\varphi_{k,j}(x)\varphi_{k,j},\quad x,y\in S^m, \quad t\in(0,\pi),
\end{eqnarray*}
we have that
\begin{eqnarray*}
S_t(K_{\sigma}(x, \cdot))(y) - K_{\sigma}(x, y) = \sum_{k=0}^{\infty}\lambda_k^{\sigma}\left(P_k^{(m-1)/2}(\cos t)-1\right)\sum_{j=1}^{d_k^m}\varphi_{k,j}(x)\varphi_{k,j}(y).
\end{eqnarray*}
An application of the usual H\"{o}lder's inequality reveals that
\begin{eqnarray*}
\left(\sum_{k=0}^{\infty}\lambda_k^{\sigma}\left(P_k^{(m-1)/2}(\cos t)-1\right)^2\sum_{j=1}^{d_k^{m}}|\varphi_{k,j}(x)|^2\right)^{1/2}\left(\sum_{k=0}^{\infty}\lambda_k^{\sigma}\sum_{j=1}^{d_k^{m}}|\varphi_{k,j}(y)|^2\right)^{1/2}.
\end{eqnarray*}
is an upper bound for $\left|S_t(K_{\sigma}(x, \cdot))(y) - K_{\sigma}(x, y)\right|$.\ The second multiplicand in the bound above is $K_{\sigma}(y,y) =1$.\ On the other hand,
since (see \cite{daidi})
$$
\left(P_k^{(m-1)/2}(\cos t)-1\right)^2\asymp \min(1,kt)^2, \quad t \in (0,\pi), \quad k \in \mathbb{Z}_+,
$$
it follows that
\begin{eqnarray*}
\left(\sum_{k=0}^{\infty}\lambda_k^{\sigma}\left(P_k^{(m-1)/2}(\cos t)-1\right)^2\sum_{j=1}^{d_k^m}|\varphi_{k,j}(x)|^2\right)^{1/2}\asymp \left(\sum_{k=0}^{\infty}\lambda_k^{\sigma}\min(1,kt)^2\sum_{j=1}^{d_k^m}|\varphi_{k,j}(x)|^2\right)^{1/2},
\end{eqnarray*}
for $t\in (0,\pi)$ and $x \in S^m$.\ Hence, we may infer that
\begin{eqnarray*}
\left|S_t(K_{\sigma}(x, \cdot))(y) - K_{\sigma}(x, y)\right| \leq C\, t^2 \left(\sum_{k=0}^{\infty} k^2 \lambda_k^{\sigma}\sum_{j=1}^{d_k^m}|\varphi_{k,j}(x)|^2\right)^{1/2},\quad t \in (0,\pi),\quad  x,y \in S^m,
\end{eqnarray*}
for some $C>0$.\ If we put
$$
B(x):= \left(\sum_{k=0}^{\infty} k^2 \lambda_k^{\sigma}\sum_{j=1}^{d_k^m}|\varphi_{k,j}(x)|^2\right)^{1/2}, \quad x\in S^m,
$$
an application of the well known addition formula for spherical harmonics yields that
$$
B(x)= \left(\sum_{k=0}^{\infty} k^2 \lambda_k^{\sigma}d_k^m\right)^{1/2}, \quad x\in S^m.
$$
Since $B$ does not depend upon $x$, in order to show that $B\in L^{\infty}(S^m)$, it suffices to verify that the series in
$$
B^2=\sigma^{m-1}e^{-2/{\sigma^2}}\Gamma((m+1)/2)\sum_{k=0}^{\infty} k^2  I_{k+(m-1)/2}(2\sigma^{-2}) d_k^m
$$
is convergent.\ However, due to the inequality
$$
I_{\nu}(x)<\frac{x^{\nu} e^x}{2^\nu \Gamma(\nu+1)}, \quad x>0,
$$
proved in \cite{luke}, we have that
$$
B^2\leq \Gamma((m+1)/2)\sum_{k=0}^{\infty} \frac{k^2 d_k^m }{\sigma^{2k}\Gamma(k+(m+1)/2)}.
$$
Introducing (\ref{ineqb}) in the expression above, it is seen that the convergence of the series boils down to the convergence of
$$
\sum_{k=k_0}^{\infty} \frac{k^{m+2}   }{\sigma^{2k}\Gamma(k+(m+1)/2)}.
$$
Basic estimates for the Gamma function reduces the analysis to the convergence of
$$
\sum_{k=k_0}^{\infty} \frac{ k^{m+2} e^{k+(m-1)/2)}  }{\sigma^{2k} [k+(m+1)/2]^{k+(m-1)/2} }.
$$
However, the series above converges by the usual ratio test.\ Thus, $B\in L^{\infty}(S^m)$ and, consequently, $K_\sigma$ is $(\mu_t, B, 2)$-H\"{o}lder with $\{\mu_t\}_{k=0}^\infty$ as in Example \ref{shifting}.
\end{ex}

\section{Approximation operators}

This section is mainly concerned with the analysis of certain normalized linear operators associated with families $\{T_t : t\in (0,\pi)\}$ of convolution operators, as defined by the constructions presented in the previous section.\ The normalized operators will be used in the search for optimal decay rates for the sequence of eigenvalues of a H\"{o}lderian integral operator $\mathcal{K}$ whenever secondary conditions are in force.\ The rates are to be used in Section 5.

Throughout this section $\{T_t: t\in(0,\pi)\}$ will denote a family of convolution operators defined by a family $\{\mu_t: t\in (0,\pi)\}$ of measures in $\mathcal{M}_{\rho}(S^m)$, in accordance with the two cases described in Section 3.\
The specific problem we will deal with here is this one: for fixed nonnegative integers $n$ and $r$ and a sequence of real or complex valued functions $\{G_n: n=1,2,\ldots\}$ which are integrable in $[0,\pi]$, to decide whether the formula
\begin{equation}\label{mainoperator}
\mathcal{A}_n(f)(\cdot):=\int_0^{\pi}G_n(t)T_t(f)(\cdot)v_m(t)(\sin t)^rdt,\quad  f \in L^p(S^m),
\end{equation}
(or a slight change of it) defines a bounded linear operator on $L^p(S^m)$.\ If it does so, to estimate the rank of such operator.\ If $T_t(f)= f * \mu_t$, $f \in L^p(S^m)$, the normalizing function $v_m(t)$ appearing above
should be interpreted as the constant function $1$ while if $T_t(f)=f* \mu_{ K_{i}^{t}}$, $f \in L^p(S^m)$, for some kernel $K^t$, then
\begin{equation}\label{vm}
v_m(t):=\int_{supp(K^t(x, \cdot))}d\sigma_m(z),
\end{equation}
where $supp(K^t(x, \cdot))$ is the support of the function $y \in S^m \to K^t(x,y)=K_i^t(x \cdot y)$.\ Observe that $supp(K^t(x, \cdot))$ does not depend upon $x$ due to the invariance of $\sigma_m$ with respect to orthogonal transformations on $\mathbb{R}^{m+1}$ (\cite{jordao}).

Particular versions of the operator given by (\ref{mainoperator}) are very commom in the approximation theory literature.\ After some acquaintance with them, we found convenient to consider (\ref{mainoperator}) in its normalized version:
\begin{equation}\label{mainoperator1}
\mathcal{A}_{n,r}(f)(\cdot):=\int_0^{\pi}G'_n(t)T_t(f)(\cdot)v_m(t)(\sin t)^rdt,  \quad f \in L^p(S^m),
\end{equation}
where $G'_n=c_{n,r}^{-1}G_n$ and
$$
c_{n,r}=\int_{0}^{\pi} |G_n(t)| v_m(t) (\sin t)^r\,dt.
$$

A sufficient condition for boundedness of $\mathcal{A}_{n,r}$ is the content of the proposition below.

\begin{prop}\label{propoperator} If $\{\mu_t: t\in (0,\pi)\}$ is uniformly bounded in $\mathcal{M}_{\rho}(S^m)$, then (\ref{mainoperator1}) defines a bounded linear operator on $L^p(S^m)$, $1\leq p\leq \infty$.\ In addition,
the family $\{\mathcal{A}_{n,r}: n=1,2,\ldots\}$ is uniformly bounded.
\end{prop}

\pf Needless to say that the linearity of (\ref{mainoperator1}) follows from that of $T_t$, $t\in (0,\pi)$.\ On the other hand, Minkowski's inequality for integrals (\cite[p.194]{folland}) implies that
$$
\|\mathcal{A}_{n,r}(f)\|_p  \leq  \int_0^{\pi}|G_n'(t)|\|T_{t}(f)\|_pv_m(t)(\sin t)^r\, dt, \quad f \in L^2(S^m).$$
Hence,
$$
\|\mathcal{A}_{n,r}(f)\|_p  \leq   \|f\|_p \int_0^{\pi}|\mu_t|  |G_n'(t)|v_m(t)(\sin t)^r\, dt\leq M \|f\|_p , \quad f \in L^p(S^m),
$$
where $M$ is a uniform upper bound for the family $\{\mu_t: t\in (0,\pi)\}$ in $\mathcal{M}_{\rho}(S^m)$.\ Thus $M$ is a uniform upper bound for the sequence $\{\|\mathcal{A}_{n,r}(f)\|_p\}_{n=1}^\infty$. \eop

A relevant realization for the operator (\ref{mainoperator1}) involves the generalized Jackson kernels $J_{l,n}$ of a fixed order $l$ defined by
$$
J_{l,n}(t):= \left[\frac{\sin((n+1)t/ 2)}{\sin(t/2)}\right]^{2l}, \quad t\in (0,\pi).
$$
In this case, $G_n'=J_{l,n}'$, where $J_{l,n}'=c_{n,r}^{-1}J_{l,n}$ and
$$c_{n,r}=\int_{0}^{\pi} \left[\frac{\sin((n+1) t/ 2)}{\sin(t/2)}\right]^{2l} v_m(t) (\sin t)^r\,dt.$$
The \textit{generalized Jackson kernel} $J_{l,n}$ is an even trigonometric polynomial of degree $ln$.\ In the case $l=1$, it reduces itself to the \textit{F\`{e}jer kernel} while the cases in which $l=1$ and $n\in 2\mathbb{Z}_+$ corresponds to the \textit{Dirichlet kernels} (\cite[p.3]{devore}).\ In all these cases, the operators $\mathcal{A}_{n,r}$ are of finite rank, and in order to prove that, we first compute the projections of $\mathcal{A}_{n,r}(f)$, for $f \in L^p(S^m)$.

\begin{prop}\label{projoperator} For $n \geq 1$ and $f$ in $L^2(S^m)$, it holds
$$
\mathcal{Y}_k(\mathcal{A}_{n,r}(f))(x)=\int_0^{\pi}G'_n(t)\mu_{t}^k \, \mathcal{Y}_k(f)(x)v_m(t)(\sin t)^rdt, \quad k=0,1, \ldots,
$$
where for each $t\in (0,\pi)$, $\{\mu_{t}^k\}_{k=0}^\infty$ is the multiplier of $\mu_t$.
\end{prop}

\pf Fix $n$ and $f$.\ It is an easy matter to verify that
$$
\mathcal{Y}_k(\mathcal{A}_{n,r}(f))(x)=\int_0^{\pi}G'_n(t)\left[\mathcal{Y}_k\left(T_t(f)\right)(x)\right]v_m(t)(\sin t)^rdt,\quad k=0,1,\ldots.
$$
On the other hand, we have that
\begin{eqnarray*}
\mathcal{Y}_k(T_t(f)) =  \mathcal{Y}_k(f\ast \mu_t)= \mu_{t}^k \, \mathcal{Y}_k(f), \quad  k=0,1,\ldots.
\end{eqnarray*}
The result follows.\eop

We close the section presenting a methodology in order to gain finite rank operators among the $\mathcal{A}_{n,r}$.\ A simplified version of the result to be described below can be found in \cite[p.214]{lizor} while another one, but in a more general setting, can be found in \cite[p.760]{platonov}.

\begin{thm}\label{rankoperator}
Let $\{\mu_t: t\in (0,\pi)\}$ be uniformly bounded in $\mathcal{M}_{\rho}(S^m)$ and assume the following assumption holds: for each $t$, the multiplier of $\mu_t$ is
$$\left\{c_{k,m}(v_m(t))^{-1}P^{\beta}_{\alpha(k)}(\cos t) (\sin t)^{\gamma}\right\}_{k=0}^\infty,$$
where $\{c_{k,m}\}_{k=0}^\infty$ is a sequence of nonzero real numbers, $\alpha: \mathbb{Z}_+ \to \mathbb{Z}_+$ is strictly increasing, $\gamma>0$ and $2\beta$ is an integer at least $\gamma$.\ If $G_n'$ is an even trigonometric polynomial of degree $\delta(n)$, then $\mathcal{A}_{n,2\beta-\gamma}$ is a bounded linear operator on $L^2(S^m)$ of rank at most $d_{\alpha^{-1}(\delta(n))}^{m+1}$.
\end{thm}
\pf We know from Proposition \ref{projoperator} that $\mathcal{A}_{n,r}$ is a bounded linear operator on $L^2(S^m)$.\ Taking into account our assumption and applying the previous proposition, we reach that
$$
\mathcal{Y}_k(\mathcal{A}_{n,r}(f))(x)=c_{k,m}\left(\int_0^{\pi}G'_n(t)P^{\beta}_{\alpha(k)}(\cos t)(\sin t)^{r+\gamma}dt\right) \mathcal{Y}_k(f)(x), \quad k=0,1, \ldots.
$$
If $G_n'$ is an even trigonometric polynomial of degree $\delta(n)$, we can write
$$
G_n'(t)=\sum_{j=0}^{\delta(n)}a_jP_j^{\beta}(\cos t),\quad a_1, a_2, \ldots, a_{\delta(n)} \in \mathbb{R}.
$$
Hence, the integrals appearing above become
$$
\int_0^{\pi}G'_n(t)P^{\beta}_{\alpha(k)}(\cos t)(\sin t)^{r+\gamma}dt = \sum_{j=0}^{\delta(n)}a_j\int_0^{\pi}P_j^{\beta}(\cos t)P^{\beta}_{\alpha(k)}(\cos t)(\sin t)^{r+\gamma}dt.
$$
We now proceed inserting the choice $r=2\beta-\gamma$, that is, we look at the integrals
$$
\int_0^{\pi}P_j^{\beta}(\cos t)P^{\beta}_{\alpha(k)}(\cos t)(\sin t)^{2\beta}dt, \quad j=0,1,\ldots, \delta(n).
$$
Since $\alpha$ is strictly increasing, we can pick $k_0$ so that $\alpha(k)> \delta(n)$ whenever $k\geq k_0$.\ Hence,
the well-known orthogonality relation for Gegenbauer polynomials (\cite[p.98]{morimoto}) implies that
$$
\int_0^{\pi}P_j^{\beta}(\cos t)P^{\beta}_{\alpha(k)}(\cos t)(\sin t)^{2\beta}dt =0, \quad \alpha(k)> \delta(n), \quad j=0,1,\ldots, \delta(n).
$$
It follows that
$$
\mathcal{Y}_k(\mathcal{A}_{n,2\beta-\gamma}(f)) = 0, \quad \alpha(k)> \delta(n),
$$
or, equivalently, that
$$
\mathcal{Y}_k(\mathcal{A}_{n,2\beta-\gamma}(f)) = 0, \quad  k>\alpha^{-1}(\delta(n)).
$$
Thus, $\mathcal{A}_{n,2\beta-\gamma}(f)$ is a polynomial of degree at most $d_{\alpha^{-1}(\delta(n))}^{m+1}$ and the proof follows. \eop

\begin{rem}Returning to Example \ref{shifting}, it is easy to see that the family of measures $\{\mu_t: t\in (0,\pi)\}$ given there fits in the setting of Theorem \ref{rankoperator} with $c_{k,m}=1$, $v_m(t)=1$ for all $t$, $\gamma=0$, $\alpha=$ the identity mapping, and $\beta=(m-1)/2$.\ Lizorkin (\cite[p. 214]{lizor}) showed that for $G_n'=J_{l,n}$, $n=1,2,\ldots$, the operator $\mathcal{A}_{n,m-1}$ has rank at most $d_{ln}^{m+1}$.

As for Example \ref{caps}, the family of measures $\{\mu_{\mathcal{Z}_{t}}: t\in (0,\pi)\}$ fits in Theorem \ref{rankoperator} with
$$
c_{k,m}=\frac{\omega_{m-1}(m-1)}{k(k+m-1)},
$$
$v_m(t)=C_m(t)$ for all $t$, $\alpha(k)=k-1$, $\beta=(m+1)/2$, $\gamma=m$.\ If $G_n'=J_{l,n}$, Proposition 2 in \cite{jordao} reveals that $\mathcal{A}_{n,1}$ has rank at most $d_{ln+1}^{m+1}$.

Regarding Example \ref{stekelov}, the family of measures $\{\mu_{\mathcal{W}_{t}}: t\in (0,\pi)\}$ also fits into the setting of Theorem \ref{rankoperator} with $c_{k,m}$ as in the previous case, $v_m(t)=D_m(t)$ for all $t$, $\alpha(k)=k-1$, $\beta=(m+1)/2$, and $\gamma=m$.\ If $G_n'=J_{l,n}$, Proposition 3 in \cite{jordao} reveals that the operator $\mathcal{A}_{n,m}$ has rank at most $d_{ln+1}^{m+1}$.
\end{rem}

\section{Kolmogorov $n$-widths}

In this section, we return to the setting and notation of Section 1.\ We will deduce decay rates for the sequence of eigenvalues of (\ref{intequ}) under the assumption that $K$ is ($\mu_t, B, \rho$)-\textit{H\"{o}lder},
in which the $\mu_t$ all belong to either setting described in Section 2.\ In a second step, we will use them to estimate $d_n(S;L^2(S^m))$.

The procedure described above is not standard.\ Indeed, usually one estimates Kolmogorov $n$-widths of certain subspaces of the reproducing kernel Hilbert space $\mathcal{H}_K$ of a smooth enough kernel $K$ in order to deduce decay rates for the sequence of eigenvalues of the integral operator $\mathcal{K}$ in a second step.\ A typical example of this standard procedure is described in \cite{santin}.

The results will be validated under an additional assumption on the normalizing function $v_m(t)$ introduced in Section 3.\ Precisely, if $T_t(f)= f * \mu_t$, $f \in L^p(S^m)$, then we will require $v_m(t)\asymp 1$, $t \in (0,\infty)$.\ Otherwise,
we will require that
 \begin{equation}\label{volumeass}
v_m(t) \asymp t^{c(m)}, \quad t\in (0,\pi),
\end{equation}
for some constante $c(m)$.

We begin with decay rates for the sequence of eigenvalues of $\mathcal{K}$, in the case in which $K$ is ($\mu_t, B, \rho$)-H\"{o}lder and some other secondary features hold.\ Prior to that, we recall
a technical lemma involving the square root $\mathcal{K}^{1/2}$ of the operator $\mathcal{K}$.\ We use $\|\cdot \|_{HS}$ to denote the usual Hilbert-Schmidt norm of a Hilbert-Schmidt operator.\ Its proof can be easily adapted
from Lemmas 4 and 5 in \cite{jordao}.

\begin{lem} \label{HSmia}
Let $\mathcal{K}$ be an integral operator as described at the introduction.\ If $K$ is ($\mu_t, B, \rho$)-H\"{o}lder, then the operator $\mathcal{K}^{1/2}- \mathcal{A}_{n,r}\circ \mathcal{K}^{1/2}$ is Hilbert-Schmidt
and there exists a positive constant $C$ so that
$$
\|\mathcal{K}^{1/2}- \mathcal{A}_{n,r}\circ \mathcal{K}^{1/2}\|_{HS}\leq C\left[\int_{0}^{\pi}G_n'(t)t^{\rho/2}v_m(t) (\sin t)^{r} \,dt\right].
$$
Further, if $\mathcal{A}_{n,r}$ has finite rank $q$, then the $(2q)$-th approximation number $a_{2q}(\mathcal{K}^{1/2})$ of $\mathcal{K}^{1/2}$
satisfies
$$
qa_{2q}(\mathcal{K}^{1/2}) \leq \|\mathcal{K}^{1/2}- \mathcal{A}_{n,r}\circ \mathcal{K}^{1/2}\|_{HS}.$$
\end{lem}

\begin{prop}\label{decaymia} Let $\{\mu_t: t\in (0,\pi)\}$ be a uniformly bounded family of measures in $M_\rho(S^m)$ so that, for each $t$, the multiplier of $\mu_t$ is
$$\left\{c_{k,m}(v_m(t))^{-1}P^{\beta}_{\alpha(k)}(\cos t) (\sin t)^{\gamma}\right\}_{k=0}^\infty,$$
where $\{c_{k,m}\}_{k=0}^\infty$ is a sequence of nonzero real numbers, $\alpha: \mathbb{Z}_+ \to \mathbb{Z}_+$ is strictly increasing, $\gamma>0$ and $2\beta$ is an integer at least $\gamma$.\
Assume that for every positive integer $l$, there exists $q=q(l,m)>0$ so that $ln\leq \alpha(qn)$, $n=1,2,\ldots$.\ If $K$ is ($\mu_t, B, \rho$)-H\"{o}lder, then the sequence of eigenvalues $\{\lambda_n\}_{n=1}^\infty$ of the integral operator $\mathcal{K}$ satisfies
$$
\lambda_n=O(n^{-1-\rho/m}), \quad n\to\infty.
$$
\end{prop}
\pf The proof begins with an application of Theorem \ref{rankoperator} with the choice
$$
G'_n=J_{l,n}', \quad n=1,2,\ldots,
$$
leaving $l$ fixed but generic.\ Due to our assumptions, we end up concluding that, for all $n$, $\mathcal{A}_{n,2\beta-\gamma}$ has rank at most $d^{m+1}_{\alpha^{-1}(ln)}$.\ Recalling the estimate for the dimensions $d_n^m$ given in Section 2, we may select a positive integer $q'$ so that
$$d^{m}_{\alpha^{-1}(ln)}\leq (q' \alpha^{-1}(ln))^m,\quad n=1,2,\ldots.$$
Invoking our assumptions on $\alpha$, we end up concluding that, for all $n$, the rank of $\mathcal{A}_{n,2\beta-\gamma}$ is at most $(qq' n)^m$.\ 
If $K$ is ($\mu_t, B, \rho$)-H\"{o}lder, we can infer from Lemma \ref{HSmia} that
$$
(qq'n)^{m}\lambda_{(qq'n)^{m}}=(qq'n)^{m}a_{(qq'n)^{m}}\leq C \left[\int_{0}^{\pi}J_{l, n}(t)t^{\rho/2}v_m(t) (\sin t)^{2\beta-\gamma} \,dt\right]^2, \quad n=1,2,\ldots,
$$
with $C>0$.This is the point where we need a special choice of the integer $l$ in order to proceed.\ Picking $l$ so that $2l \geq \rho+c(m)+2\beta-\gamma+1$ ($c(m)$ is the constant in (\ref{volumeass})), we can apply Lemma 1 in \cite{jordao} in order to see that
$$
(qq'n)^{m}\lambda_{(qq'n)^{m}}\leq \frac{l^{\rho}\, C_{m,\rho}^{l} }{n^{\rho}},\quad n \in l\mathbb{Z}_+,
$$
in which $C_{m,\rho}^{l}$ is now a positive constant depending upon $m,\rho,l$.\
Going one step further, we may repeat the trick for $n \in j+l\mathbb{Z}_+$, $j \in \{1,2,\ldots, l-1\}$, and finally deduce that
$$
(qq'n)^{m}\lambda_{(qq'n)^{m}}\leq \frac{C'}{n^{\rho}}, \quad n=1,2,\ldots.
$$
for some positive constant $C'$ depending on $m,\rho,l$, and $c(m)$, but not on $n$.\ In other words,
$$
\lambda_{(qq'n)^{m}}\leq \frac{C''}{n^{\rho +m}},\quad n=1,2,\ldots.
$$
for some positive constant $C''$ not depending upon $n$.\ This implies the eigenvalue behavior described in the statement of the theorem.\eop

The extra assumption we have made on the mapping $\alpha$ in the previous theorem is not unreal.\ Indeed, it is obviously true if $\alpha$ is a affine mapping, a fact in the concrete examples quoted so far.

Taking into account Proposition \ref{decaymia} and one the formulas quoted at the introduction, we have our final result of the section.

\begin{thm}\label{Kolmogorovn}
Under the assumptions in Proposition \ref{decaymia}, if $K$ is ($\mu_t, B, \rho$)-H\"{o}lder, then
$$
d_n(S;L^2(S^m))= O((n+1)^{-1/2-\rho/2m}), \quad n\to\infty.
$$
\end{thm}

We now return to some of the examples we previously mentioned.\ If the family $\{\mu_t \in (0\,\pi)\}$ is given via the shifting operator as in Example \ref{shifting}, then due to all the comments we have made along the text, it is easily seen that the previous theorems hold for a
($\mu_t, B, \rho$)-H\"{o}lder kernel $K$.\ The same is true for the sequence $\{\mu_{\mathcal{Z}_{t}}: t \in (0,\pi)\}$ in Example \ref{caps} attached to the average on caps operators.\ Indeed, in this case, we need to observe that
$$
v_m(t)\asymp t^m, \quad t\in (0,\pi),
$$
a consequence of the inequality
$$
\frac{\tau_{m-1}}{m}\left(\frac{2}{\pi}\right)^{m-1} t^m\leq v_m(t) \leq \tau_{m-1}t^m,  \quad t \in (0,\pi),
$$
and some calculations (\cite[Example 1]{jordao}).\ Then the average on caps also fits into the assumption made in formula (\ref{volumeass}) with $c(m)=m$, $m\in \mathbb{Z}_+$.\ Finally, since the measures in Example \ref{stekelov} are related to those in  Example \ref{caps}, the same is true for them.\ Thus, in all three cases, the asymptotic behavior in Theorem \ref{Kolmogorovn} holds.

Let us finish the paper with an unusual example.\ Here we will assume $m\geq 2$ and will consider the dot product kernel
\begin{equation}\label{dotkernel}
K(x,y)=1+\sum_{k=1}^{\infty}\left(\frac{2^{k+1}n^{(m-1)/2}}{k^{1+k\epsilon/m}}\right)(x\cdot y)^n, \quad x,y\in S^m,
\end{equation}
where $\epsilon$ is chosen to be strictly bigger than $m/2$.\ If
$$
b_k:=\frac{2^{k+1}k^{(m-1)/2}}{k^{1+k\epsilon/m}}, \quad k=0,1,\ldots,
$$
then it is easily seen that
$$
k^{2\epsilon}\frac{b_k}{b_{k-1}} \to 2e^{-\epsilon/m},\quad  k\to\infty.
$$
In other words, $K$ satisfies all the assumptions of Theorem 3.3 in \cite{azevedo}.\ In particular, the sequence $\{\lambda_k\}_{k=1}^\infty$ of eigenvalues of $\mathcal{K}$ satisfies
$$
\lambda_k \asymp\frac{b_k}{2^{k+1}k^{(m-1)/2}} = \frac{1}{k^{1+ k\epsilon/m}},\quad   k\to\infty.
$$
Taking into account Example \ref{gaussianex} and considering the Mercer expansion of $K$, it can be seen that there exists a positive constant $c$, such that
\begin{eqnarray*}
\left|S_t(K(x, \cdot))(y) - K(x, y)\right| \leq c\, t^2 \left(\sum_{k=0}^{\infty} k^2 \lambda_k\sum_{j=1}^{d_k^m}|\varphi_{k,j}(x)|^2\right)^{1/2},\quad x,y \in S^m.
\end{eqnarray*}
Defining
$$
B(x):= \left(\sum_{k=0}^{\infty} k^2 \lambda_k\sum_{j=1}^{d_k^m}|\varphi_{k,j}(x)|^2\right)^{1/2}, \quad x\in S^m,
$$
we immediately have that
$$
B(x)\leq \left(\sum_{k=0}^{\infty} \lambda_k(k \, d_k^m)^2\right)^{1/2} \asymp \left(\sum_{k=0}^{\infty} \frac{(k \, d_k^m)^2}{k^{1+ k\epsilon/m}}\right)^{1/2}, \quad x\in S^m.
$$
Since the series appearing above is clearly convergent we conclude that $B\in L^{\infty}(S^m)$.\
It is now clear that $K$ is $(S_t, B,2)$-H\"{o}lder and, as so, Theorem \ref{Kolmogorovn} is applicable, the outcome being
$$
d_n(S;L^2(S^m))  \asymp \frac{1}{(n+1)^{1/2+ (n+1)\epsilon/2m}},\quad   n\to\infty.
$$
In particular, it follows that
$$
d_n(S;L^2(S^m)) = o\left((n+1)^{-\left(1/2+(n+1)/2m\right)}\right),\quad   n\to\infty.
$$
%
%---------------------------------------------------------------------------------------------------------------------------------------------
%

%
%----------------------------------------------------------------------------------------------------------------------------------------------------------
%

\vspace*{1.5cm}

\noindent %T. Jord\~{a}o and V. A. Menegatto \\
Departamento de
Matem\'atica,\\ ICMC-USP - S\~ao Carlos, Caixa Postal 668,\\
13560-970 S\~ao Carlos SP, Brasil\\ E-mails: tjordao@icmc.usp.br; menegatt@icmc.usp.br

\end{document}